# Partial Zeta Functions of Algebraic Varieties over Finite Fields

Daqing Wan

Department of Mathematics

University of California

Irvine, CA 92697-3875

dwan@math.uci.edu

*Dedicated to Professor Chao Ko for his 90-th birthday*

**Abstract**. Motivated by arithmetic applications, we introduce the notion of a partial zeta function which generalizes the classical zeta function of an algebraic variety defined over a finite field. We then explain two approaches to the general structural properties of the partial zeta function in the direction of the Weil type conjectures. The first approach, using an inductive fibred variety point of view, shows that the partial zeta function is rational in an interesting case, generalizing Dwork's rationality theorem. The second approach, due to Faltings, shows that the partial zeta function is always nearly rational.

## 1. Introduction

Let $\mathbf{F}_q$ be the finite field of $q$ elements of characteristic $p$. Let $X$ be an affine algebraic variety over $\mathbf{F}_q$, embedded in some affine space $\mathbf{A}^n$. That is, $X$ is defined by a system of polynomial equations

$$F_1(x_1, \cdots, x_n) = \cdots = F_r(x_1, \cdots, x_n) = 0,$$

where each $F_i$ is a polynomial defined over $\mathbf{F}_q$. Let $d_1, \cdots, d_n$ be $n$ positive integers. For each positive integer $k$, let

$$X_{d_1, \cdots, d_n}(k, X) = \{x \in X(\bar{\mathbf{F}}_q) | x_1 \in \mathbf{F}_{q^{d_1 k}}, \cdots, x_n \in \mathbf{F}_{q^{d_n k}}\}.$$

That is, the degree over $\mathbf{F}_q$ of the $i$-th coordinate $x_i$ of $x = (x_1, \cdots, x_n) \in X$ divides $d_i k$ for $1 \leq i \leq n$. Let

$$N_{d_1, \cdots, d_n}(k, X) = \# X_{d_1, \cdots, d_n}(k, X).$$

The number $N_{d_1, \cdots, d_n}(k, X)$ counts the points of $X$ whose coordinates lie in different subfields of $\bar{\mathbf{F}}_q$. We would like to understand this sequence of integers $N_{d_1, \cdots, d_n}(k, X)$ as $k$ varies. As usual, it is sufficient to understand the following generating function.

**Definition 1.1.** *Given $n$ positive integers $d_1, \cdots, d_n$. The associated partial zeta function $Z_{d_1, \cdots, d_n}(X, T)$ of $X/\mathbf{F}_q$ is defined to be the following formal power series*

$$Z_{d_1, \cdots, d_n}(X, T) = \exp(\sum_{k=1}^{\infty} \frac{N_{d_1, \cdots, d_n}(k, X)}{k} T^k).$$

In the special case that $d_1 = \cdots = d_n = d$, the number $N_{d, \cdots, d}(k, X)$ is just the number of $\mathbf{F}_{q^{dk}}$-rational points on $X$. The partial zeta function $Z_{d, \cdots, d}(X, T)$ then becomes the classical zeta function $Z(X \otimes \mathbf{F}_{q^d}, T)$ of the variety $X \otimes \mathbf{F}_{q^d}$, which is a rational function by Dwork's rationality theorem [Dw1] and satisfies a suitable Riemann hypothesis by Deligne's theorem [De2] on the Weil conjectures.

Our motivation to introduce the above more general partial zeta function comes directly from potential applications in number theory, combinatorics and coding theory. The idea of counting only part of the rational points corresponds to the classical partial character sum problem in number theory. For arbitrary partial counting, the problem would be too difficult and one cannot expect a good mathematical structure for it. Our partial counting is not arbitrary as each coordinate $x_i$ is restricted to run over a subfield (not an arbitrary subset of $\bar{\mathbf{F}}_q$). Thus, it seems reasonable to expect good mathematical properties of the partial zeta function $Z_{d_1,\cdots,d_n}(X,T)$, generalizing what we know about the Weil conjectures for the classical zeta function $Z(X,T)$. For certain plane curves of Artin-Schreier and Kummer type, it was observed in [W1] that the partial counting problem (as well as certain partial one variable character sum problem) has a fairly satisfactory solution by exploiting Weil's work on the Riemann hypothesis for curves over finite fields. In the present paper, we consider the general and higher dimensional case using the theory of $\ell$-adic cohomology. The situation is naturally more complicated. We have

**Theorem 1.2.** *Let $\{d_1, d_2, \cdots, d_n\}$ be a sequence of positive integers which can be rearranged such that $d_1|d_2|\cdots|d_n$. Then the partial zeta function $Z_{d_1,\cdots,d_n}(X,T)$ is a rational function satisfying a suitable Riemann hypothesis. That is, there are finitely many algebraic integers $\alpha_a$ and $\beta_b$ depending on $X/\mathbf{F}_q$ and $\{d_1,\cdots,d_n\}$ such that*

$$N_{d_1,\cdots,d_n}(k,X) = \sum_a \alpha_a^k - \sum_b \beta_b^k$$

*for every positive integer $k$, where the absolute value of each $\alpha_a$ (resp. $\beta_b$) is an integral power of $\sqrt{q}$.*

This result came to my attention in late 93 in connection with my initial study of Dwork's conjecture [Dw2] on the $p$-adic meromorphic continuation of his unit root zeta function. The idea is to view $X$ as a sequence of fibred varieties

$$X = X_n \xrightarrow{f_n} X_{n-1} \xrightarrow{f_{n-1}} X_{n-2} \cdots \xrightarrow{f_2} X_1$$

such that

$$N_{d_1,\cdots,d_n}(k,X) = \#\{x \in X_n(\mathbf{F}_{q^{d_n k}}) | f_{j+1} \circ \cdots \circ f_n(x) \in X_j(\mathbf{F}_{q^{d_j k}}), 1 \le j \le n\}.$$

The desired rationality of the partial zeta function then follows from an inductive (or iterated) argument using Newton's formula and the theory of $\ell$-adic cohomology. From fibred variety point of view, the varieties $X_i$ do not have to be affine. They can be any scheme of finite type over $\mathbf{F}_q$. Our proof works for more general partial L-functions attached to a character.

If the integers $d_i$ do not satisfy the dividing condition, the partial zeta function should not be rational in general but should be nearly rational in some sense. I explained this problem to Gerd Faltings in the summer of 99 when I was visiting the Max-Planck Institute. Faltings showed me a clever geometric construction which reduces the problem to a "more general" Lefschtez fixed point theorem in $\ell$-adic cohomology. His idea is to construct a new variety $Y(d_1,\cdots,d_n,X)/\mathbf{F}_q$ with an automorphism $\sigma$ of order $d = [d_1,\cdots,d_n]$ (the least common multiple) such that

$$N_{d_1,\cdots,d_n}(k,X) = \#\{y \in Y(d_1,\cdots,d_n,X)(\bar{\mathbf{F}}_q) | \sigma \circ \mathrm{Frob}^k(y) = y\},$$

where Frob denotes the $q$-th power Frobenius map acting on $Y$. The general $\ell$-adic trace formula then implies that the partial zeta function is indeed nearly rational in the following precise sense.

**Theorem 1.3 (Faltings).** *Let $\{d_1,\cdots,d_n\}$ be $n$ positive integers. Let $d = [d_1,\cdots,d_n]$ be the least common multiple of the $d_i$. Let $\zeta_d$ be a primitive $d$-th root of unity. There are $d$ rational functions $R_j(T)$ $(1 \le j \le d)$ with $R_j(0) = 1$ and with algebraic integer coefficients such that*

$$Z_{d_1,\cdots,d_n}(X,T) = \prod_{j=1}^d R_j(T)^{\zeta_d^j}.$$

Furthermore, each rational function $R_j(T)$ satisfies a suitable Riemann hypothesis. That is, there are finitely many algebraic integers $\alpha_{aj}$ and $\beta_{bj}$ depending on $X/\mathbf{F}_q$ and $\{d_1, \cdots, d_n\}$ such that

$$N_{d_1, \cdots, d_n}(k, X) = \sum_{j=1}^{d} \zeta_d^j (\sum_a \alpha_{aj}^k - \sum_b \beta_{bj}^k)$$

for every positive integer $k$, where the absolute value of each $\alpha_{aj}$ (resp. $\beta_{bj}$) is an integral power of $\sqrt{q}$.

This theorem shows that for many purposes, the partial zeta function is as good as a rational function. In the special case that $d \leq 2$, then $\zeta_d^j = \pm 1$ and the partial zeta function is actually rational, which also follows from Theorem 1.2. The author does not know if one can prove Theorem 1.2 by Faltings' method. It should be noted that we do not claim that the rational functions $R_j(T)$ in Theorem 1.3 have rational integer coefficients, which would imply that the partial zeta function is rational. However, presumably, the rational functions $R_j(T)$ have coefficients which are integers in the $d$-th cyclotomic field $\mathbf{Q}(\zeta_d)$. Although the partial zeta function is not expected to be a rational function in general, the author has not found a counter-example.

The above results show that there are powerful theoretical methods available to handle partial zeta functions. Virtually, all theories developed for the Weil conjectures, $p$-adic or $\ell$-adic (particularly Deligne's main theorem), can be used to study the partial zeta function. For practical applications, however, one needs to have precise and sharp information on the weights and slopes about the zeros and poles of the partial zeta function. This can be quite difficult and complicated in general, going beyond the Weil conjectures even for some geometrically very simple varieties. The new arithmetic parameters $\{d_1, \cdots, d_n\}$ allow one to ask many questions about the variation of the partial zeta function $Z_{d_1, \cdots, d_n}(X, T)$ when some of the arithmetic parameters $\{d_1, \cdots, d_n\}$ vary. These variation questions seems to have important arithmetic meaning, see section 4 for further remarks and questions.

As a preliminary $p$-adic result, we remark that one can easily obtain the following theorem from the Ax-Katz theorem [Ka]. The proof is left to the reader as an excercise.

**Theorem 1.4.** *Let $d = [d_1, \cdots, d_n]$. Let $D_j$ $(1 \leq j \leq r)$ be the degree of the polynomial $F_j$. Let $\mu(d_1, \cdots, d_n; D_1, \cdots, D_r)$ denote the smallest non-negative integer which is at least as large as*

$$\frac{d_1 + \cdots + d_n - (D_1 + \cdots + D_r)d}{\max_{1 \leq j \leq r} D_j}.$$

*Then for every positive integer $k$, we have the $p$-adic estimate*

$$\mathrm{ord}_q N_{d_1, \cdots, d_n}(k, X) \geq k\mu(d_1, \cdots, d_n; D_1, \cdots, D_r).$$

*Equivalently, for every reciprocal zero (resp. reciprocal pole) $\alpha$ occuring in Theorems 1.2-1.3, we have the $p$-adic estimate*

$$\mathrm{ord}_q(\alpha) \geq \mu(d_1, \cdots, d_n; D_1, \cdots, D_r).$$

mailx

## 2. Rationality of partial zeta functions

In this section, we prove Theorem 1.2 and its generalization to partial L-functions attached to a constructible $\ell$-adic sheaf, where $\ell$ is a prime number different from $p$. The word "variety" over $\mathbf{F}_q$ in this paper means a separated scheme of finite type over $\mathbf{F}_q$.

For a variety $X$ over $\mathbf{F}_q$ and a positive integer $k$, we define

$$X(k) = X(\mathbf{F}_{q^k})$$

to be the set of $\mathbf{F}_{q^k}$-rational points on $X$. Let

$$X = X_n \xrightarrow{f_n} X_{n-1} \xrightarrow{f_{n-1}} X_{n-2} \cdots \xrightarrow{f_2} X_1$$

be a sequence of fibred varieties over $\mathbf{F}_q$. That is, the varieties $X_i$ and the morphisms $f_i$ are all defined over $\mathbf{F}_q$. For $1 \leq i < j \leq n$, define

$$g_{ij} = f_{i+1} \circ \cdots \circ f_j : X_j \longrightarrow X_i.$$

Let $d_1 | d_2 | \cdots | d_n$ be a sequence of positive integers such that each dividing the next. For $1 \leq j \leq n$, define

$$X_{d_1, \cdots, d_j}(k, X_j) = \{x \in X_j(d_j k) | g_{ij}(x) \in X_i(d_i k), 1 \leq i < j\}.$$

Thus, for $j > 1$, we have

$$X_{d_1, \cdots, d_j}(k, X_j) = \{x \in X_j(d_j k) | f_j(x) \in X_{d_1, \cdots, d_{j-1}}(k, X_{j-1}).$$

Let $\mathcal{E}_j$ be a constructible $\ell$-adic sheaf on $X_j / \mathbf{F}_q$. Let Frob denote the Frobenius map on $\mathcal{E}_j$ induced by the $q$-th power Frobenius map on $X_j$. For positive integers $k$ and $h$, we define the partial character sum

$$S_{d_1, \cdots, d_j}(\mathcal{E}_j^h, X_j(d_j k)) = \sum_{x \in X_{d_1, \cdots, d_j}(k, X_j)} \mathrm{Tr}(\mathrm{Frob}^{h d_j k} | \mathcal{E}_{jx}), \tag{2.1}$$

where $\mathcal{E}_{jx}$ denotes the fibre of $\mathcal{E}_j$ at the geometric point $x$ of $X_j$. In particular, taking $j = 1$, we have

$$S_{d_1}(\mathcal{E}_1^h, X_1(d_1 k)) = \sum_{x \in X_1(d_1 k)} \mathrm{Tr}(\mathrm{Frob}^{h d_1 k} | \mathcal{E}_{1x}).$$

Taking $h = 1$, the above equation becomes the following standard complete character sum over $X_1(d_1 k)$:

$$S_{d_1}(\mathcal{E}_1, X_1(d_1 k)) = \sum_{x \in X_1(d_1 k)} \mathrm{Tr}(\mathrm{Frob}^{d_1 k} | \mathcal{E}_{1x}). \tag{2.2}$$

The L-function attached to the sequence of character sums in (2.2) is rational by Grothendieck's rationality theorem [Gr], where the base variety is $X_1 \otimes \mathbf{F}_{q^{d_1}}$. We want to extend this result to the L-function attached to the sequence of sums in (2.1).

We have

**Theorem 2.1.** *Let $d_1 | d_2 | \cdots | d_n$. Let $\mathcal{E}$ be a constructible $\ell$-adic sheaf on $X_n$. Let $h$ be another positive integer. There are a finite number of integers $a_m$ and constructible $\ell$-adic sheaves $\mathcal{F}_m$ on $X_1$ such that for every positive integer $k$, we have the formula*

$$S_{d_1, \cdots, d_n}(\mathcal{E}^h, X_n(d_n k)) = \sum_m a_m S_{d_1}(\mathcal{F}_m, X_1(d_1 k)). \tag{2.3}$$

*Furthermore, if $\mathcal{E}$ is mixed on $X_n$ of weight at most $w$, then each $\mathcal{F}_m$ is mixed on $X_1$ of weight at most*

$$h d_n w + 2 \sum_{i=2}^{n} d_i (\dim X_i - \dim X_{i-1}).$$

**Proof**. We shall use the following universal formula

$$\text{Tr}(\phi^h|V) = \sum_{s=1}^{h} (-1)^{s-1} s \cdot \text{Tr}(\phi|\text{Sym}^s V \otimes \wedge^{h-s} V),$$  (2.4)

which holds for any linear map $\phi$ acting on a finite dimensional vector space $V$. This formula can be easily proved, see Lemma 4.1 in [W2]. It is an improvement of Newton's formula on symmetric functions. The universal formula in (2.4) shows that without loss of generality, we can assume that $h = 1$ in proving Theorem 2.1.

Now, with $h = 1$, we can write

$$S_{d_1,\cdots,d_n}(\mathcal{E}, X_n(d_n k)) = \sum_{x \in X_{d_1,\cdots,d_{n-1}}(k, X_{n-1})} \sum_{y \in f_n^{-1}(x)(d_n k)} \text{Tr}(\text{Frob}^{d_n k}|\mathcal{E}_y),$$

where the fibre $f_n^{-1}(x)$ is viewed as a variety defined over $\mathbf{F}_{q^{d_{n-1}k}}$ and $f_n^{-1}(x)(d_n k)$ denotes the set of $\mathbf{F}_{q^{d_n k}}$-rational points on $f_n^{-1}(x)$. Applying the relative $\ell$-adic trace formula to the morphism $f_n : X_n \to X_{n-1}$, we deduce that for all $x \in X_{d_1,\cdots,d_{n-1}}(k, X_{n-1})$,

$$\sum_{y \in f_n^{-1}(x)(d_n k)} \text{Tr}(\text{Frob}^{d_n k}|\mathcal{E}_y) = \sum_{m \geq 0} (-1)^{m-1} \text{Tr}(\text{Frob}^{d_{n-1}k \frac{d_n}{d_{n-1}}}|(R^m f_{n!}\mathcal{E})_x),$$

where $R^m f_{n!}\mathcal{E}$ denotes the $\ell$-adic higher direct image with compact support, which is a constructible $\ell$-adic sheaf on $X_{n-1}$. It follows that

$$S_{d_1,\cdots,d_n}(\mathcal{E}, X_n(d_n k)) = \sum_{m \geq 0} (-1)^{m-1} S_{d_1,\cdots,d_{n-1}}((R^m f_{n!}\mathcal{E})^{\frac{d_n}{d_{n-1}}}, X_{n-1}(d_{n-1}k)).$$

Iterating this procedure, or by induction, we get the required constructible sheaves $\mathcal{F}_m$ and the integers $a_m$ such that (2.3) holds. If $\mathcal{E}$ is mixed of weight at most $w$, Deligne's main theorem [De2] shows that each higher direct image sheaf $R^m f_{n!}\mathcal{E}$ of $\mathcal{E}$ is mixed of weight at most $w + 2(\dim X_n - \dim X_{n-1})$. This together with our inductive proof shows that each $\mathcal{F}_m$ is also mixed. The upper bound on the weights of $\mathcal{F}_m$ follows from the above inductive proof. It can also be checked directly from a trivial archimedian estimate. The proof of Theorem 2.1 is complete.

**Theorem 2.2.** *With the same assumption as in Theorem 2.1, the partial L-function*

$$L_{d_1,\cdots,d_n}(\mathcal{E}^h, T) = \exp\left(\sum_{k=1}^{\infty} \frac{S_{d_1,\cdots,d_n}(\mathcal{E}^h, X_n(d_n k))}{k} T^k\right)$$

*is a rational function. If $\mathcal{E}$ is mixed of integral weights, then the partial L-function $L_{d_1,\cdots,d_n}(\mathcal{E}^h, T)$ satisfies a suitable Riemann hypothesis.*

**Proof**. Applying the usual $\ell$-adic trace formula to the right side of (2.3), we obtain

$$L_{d_1,\cdots,d_n}(\mathcal{E}^h, T) = \prod_m L(\mathcal{F}_m/X_1 \otimes \mathbf{F}_{q^{d_1}}, T)^{a_m}$$

By Grothendieck's rationality theorem, each classical L-function $L(\mathcal{F}_m/X_1 \otimes \mathbf{F}_{q^{d_1}}, T)$ is a rational function. Since the powers $a_m$ are integers, we conclude that the partial L-function $L_{d_1,\cdots,d_n}(\mathcal{E}^h, T)$ is indeed rational.

Theorem 1.2 is the special case of Theorem 2.2 by taking $\mathcal{E}$ to be the constant sheaf $\mathbf{Z}_\ell$, $X_n = X$ and $X_i = \mathbf{A}^i$ for $1 \leq i < n$, with the map $f_j$ $(1 < j \leq n)$ being the natural projection

$$f_j : (x_1, \cdots, x_j) \in X_j \longrightarrow (x_1, \cdots, x_{j-1}) \in X_{j-1}.$$

In order to get a good estimate on the weights of the zeros and poles of the partial zeta function in various cases, the above inductive proof shows that we need to understand how geometrically constant sheaves occur in various combinations of tensor products and higher direct images arising from $\mathcal{E}$. This is related to representation theory of the various monodromy groups of the sheaves arising from the above inductive proof.

**Remark**. Deligne's conjecture [De2] says that every constructible sheaf is mixed. This seems still open.

## 3. Near rationality of partial zeta functions

In this section, we explain Faltings' proof of his near rationality theorem (Theorem 1.3). The notations are the same as in the introduction section.

Let $d = [d_1, \cdots, d_n]$ be the least common multiple of the positive integers $d_1, \cdots, d_n$. Let $X$ be an affine variety over $\mathbf{F}_q$, embedded in $\mathbf{A}^n$. The $d$-fold product $X^d$ of $X$ has two actions. One is the $q$-th power Frobenius action denoted by Frob. Another is the automorphism $\sigma$ on $X^d$ defined by the cyclic shift

$$\sigma(y_1, \cdots, y_d) = (y_d, y_1, \cdots, y_{d-1}),$$

where $y_j$ denotes the $j$-th component $(1 \leq j \leq d)$ of a point $y = (y_1, \cdots, y_d)$ on the $d$-fold product $X^d$. Thus, each component $y_j$ is a point on $X$, not the $j$-th coordinate in $\mathbf{A}^n$ of some point on $X$. Write

$$y_j = (y_{1j}, \cdots, y_{nj}),$$

where $y_{ij}$ $(1 \leq i \leq n)$ is the $i$-th coordinate in $\mathbf{A}^n$ of the point $y_j \in X$. For each $1 \leq i \leq n$, define a morphism

$$z_i : X^d \longrightarrow \mathbf{A}^d$$

by

$$z_i(y) = z_i(y_1, \cdots, y_d) = (y_{i1}, y_{i2}, \cdots, y_{id}).$$

Let $Y = Y(d_1, \cdots, d_n, X)$ be the subvariety of $X^d$ defined by

$$z_i \circ \sigma^{d_i} = z_i$$

for all $1 \leq i \leq n$. Thus, a point $y = (y_1, \cdots, y_d) \in X^d$ is on the subvariety $Y(d_1, \cdots, d_n, X)$ if and only if

$$y_{ij} = y_{i(j+d_i)}, \ 1 \leq i \leq n, \ 1 \leq j \leq d, \tag{3.1}$$

where $j + d_i$ is taken to be the smallest positive residue of $j + d_i$ modulo $d$. It is clear that $Y$ is stable under the action of $\sigma$.

Let $y = (y_1, \cdots, y_d)$ be a geometric point of $Y(d_1, \cdots, d_n, X)$. One checks that

$$\sigma \circ \mathrm{Frob}^k(y) = y \Longleftrightarrow y_j^{q^k} = y_{j+1}, \ 1 \leq j \leq d, \tag{3.2}$$

equivalently, if and only if

$$y_{ij}^{q^k} = y_{i(j+1)}, \ 1 \le i \le n, \ 1 \le j \le d. \tag{3.3}$$

Iterating equation (3.3) $d_i$ times, we get

$$y_{ij}^{q^{d_i k}} = y_{i(j+d_i)}.$$

Since $y$ is on $Y(d_1, \cdots, d_n, X)$, by (3.1), we deduce that

$$y_{ij}^{q^{d_i k}} = y_{ij}.$$

This shows that every fixed point $y$ of $\sigma \circ \mathrm{Frob}^k$ acting on $Y(d_1, \cdots, d_n, X)$ is uniquely determined by $y_1 \in X$ such that the $i$-th coordinate $y_{i1}$ of $y_1 \in X$ is in the field $\mathbf{F}_{q^{d_i k}}$ for all $1 \le i \le n$. That is,

$$\sigma \circ \mathrm{Frob}^k(y) = y \Longleftrightarrow y_1 \in X_{d_1, \cdots, d_n}(k, X).$$

Thus,

$$N_{d_1, \cdots, d_n}(k, X) = \#\{y \in Y(d_1, \cdots, d_n, X)(\bar{\mathbf{F}}_q) | \sigma \circ \mathrm{Frob}^k(y) = y\}.$$

The $\ell$-adic trace formula implies that

$$N_{d_1, \cdots, d_n}(k, X) = \sum_{a=0}^{2\dim(Y)} (-1)^{a-1} \mathrm{Tr}(\sigma \circ \mathrm{Frob}^k | H_c^a(Y \otimes \bar{\mathbf{F}}_q),$$

where $H_c^a$ denotes the $\ell$-adic cohomology with compact support. As $\sigma$ commutes with Frob, we can write

$$\mathrm{Tr}(\sigma \circ \mathrm{Frob}^k | H_c^a(Y \otimes \bar{\mathbf{F}}_q) = \sum_j \alpha_{aj} \beta_{aj}^k,$$

where each $\alpha_{aj}$ is an eigenvalue of $\sigma$ and each $\beta_{aj}$ is an eigenvalue of Frob. Since $\sigma^d = 1$, each $\alpha_{aj}$ is a $d$-th root of unity. By Deligne's theorem [De2], each $\beta_{aj}$ has complex absolute value $q^{w_{aj}/2}$, where $w_{aj}$ is an integer satisfying $0 \le w_{aj} \le a \le 2\dim(Y)$. Theorem 1.3 is proved.

To get further and more precise results in various interesting special cases, one needs to understand the geometry of the affine variety $Y(d_1, \cdots, d_n, X)$, such as its dimension, geometric irreducibility, number of components, etc. None of these simple questions seems to have a clean answer. Since $Y$ is affine, we do know, by the Lefschetz affine theorem, that $H_c^a(Y \otimes \bar{\mathbf{F}}_q) = 0$ for $0 \le a < \dim(Y)$. Ideally, in certain nice cases, one would hope that $H_c^a(Y \otimes \bar{\mathbf{F}}_q) = 0$ for $\dim(Y) < a < 2\dim(Y)$.

In the case that $d_1 | d_2 | \cdots | d_n$, Theorem 1.2 suggests that the eigenvalues of $\sigma$ acting on $H_c^a(Y \otimes \bar{\mathbf{F}}_q)$ might be $\pm 1$. Is this true? If so, it would give a new proof of Theorem 1.2 and perhaps some new information as well. If not, it would give some interesting cancellation information on the Frobenius eigenvalues. Theorem 1.3 can be generalized somewhat. We state one generalization here. The proof is left to the reader since it is completely similar to the above proof.

Let $\{d_1, \cdots, d_n\}$ be $n$ positive integers. For each $1 \le i \le n$, let

$$f_i : X \longrightarrow X_i$$

be a morphism of algebraic varieties defined over $\mathbf{F}_q$. For each positive integer $k$, let

$$X_{d_1, \cdots, d_n}(k, f) = \{x \in X(\bar{\mathbf{F}}_q) | f_i(x) \in X_i(\mathbf{F}_{q^{d_i k}}), 1 \le i \le n\}.$$

Let

$$N_{d_1, \cdots, d_n}(k, f) = \# X_{d_1, \cdots, d_n}(k, f).$$

Assume that this number is finite for all $k$, which is the case if and only if the map

$$f : X \longrightarrow X_1 \times \cdots \times X_n, \ f(x) = (f_1(x), \cdots, f_n(x))$$

is set theoretically finite. Define

$$Z_{d_1, \cdots, d_n}(f, T) = \exp(\sum_{k=1}^{\infty} \frac{N_{d_1, \cdots, d_n}(k, f)}{k} T^k).$$

**Theorem 3.1.** *Assume that the map*

$$f : X \longrightarrow X_1 \times \cdots X_n, \ f(x) = (f_1(x), \cdots, f_n(x))$$

*is set theoretically injective. Then, there are $d$ rational functions $R_j(T)$ $(1 \leq j \leq d)$ with $R_j(0) = 1$ and with algebraic integer coefficients such that*

$$Z_{d_1, \cdots, d_n}(f, T) = \prod_{j=1}^{d} R_j(T)^{\zeta_d^j}.$$

*Furthermore, each rational function $R_j(T)$ satisfies a suitable Riemann hypothesis.*

Theorem 1.3 is the special case of Theorem 3.1 by taking $X \hookrightarrow \mathbf{A}^n$, $X_i = \mathbf{A}^1$ and the map $f_i$:

$$f_i : (x_1, \cdots, x_n) \in X \longrightarrow x_i \in X_i = \mathbf{A}^1.$$

In this case, it is clear that the map

$$f : (x_1, \cdots, x_n) \in X \longrightarrow (x_1, \cdots, x_n) \in \mathbf{A}^n$$

is an embedding. It may be of interest to know if Theorem 3.1 can be further generalized by assuming only that the map $f$ is finite.

## 4. Further Remarks

In this final section, we briefly discuss some examples and their connections with other more classical problems.

For a concrete example, let $X$ be the rational surface

$$x_1^2 = x_2(x_2 - 1)(x_2 - x_3)$$

in $\mathbf{A}^3$, which is also the universal level 2 elliptic curve parametrized by $x_3$. In this case, the partial zeta function $Z_{d,d,1}(X, T)$ is rational by Theorem 1.2. One can show that $Z_{d,d,1}(X, T)$ is essentially the quotient of two Hecke polynomials. The integer $d$ in this example corresponds to the weight of modular forms. Although the degree of each Hecke polynomial is known in terms of dimensions of modular forms, the exact number of zeros (resp. poles) of $Z_{d,d,1}(X, T)$ is not known due to the possible cancellation of zeros. This seems to be a deep problem. The precise weight estimate of zeros and poles of $Z_{d,d,1}(X, T)$ is equivalent to the Ramanujan-Peterson conjecture. Very little is known about the slopes of the zeros and poles of $Z_{d,d,1}(X,T)$. The slope estimate of the zeros and poles of $Z_{d,d,1}(X,T)$ is closely related to the $p$-adic Ramanujan-Peterson conjecture. The $p$-adic variation of the slopes as $d$ varies $p$-adically is closely related to Gouvêa-Mazur's conjecture on modular forms, see [W2] for further information on these $p$-adic questions.

For a more general but still concrete example, let $X$ be an affine hypersurface in $\mathbf{A}^n$. Let $d_1|d_2\cdots|d_n$. The partial zeta function $Z_{d_1,\cdots,d_n}(X,T)$ is rational by Theorem 1.2. A simple heuristic argument shows that on the average, one would expect that

$$N_{d_1,\cdots,d_n}(k,X) = q^{k(d_1+\cdots+d_{n-1})} + O(q^{k(d_1+\cdots+d_{n-1})/2}). \qquad (4.1)$$

It would be interesting to establish such sharp estimates in various important cases arising from applications. For instance, if $d_i = d$ for all $1 \leq i \leq n$, then (4.1) is true for a sufficiently smooth affine hypersurface $X$. This is a consequence of Deligne's theorem [De1]. If $d_1 = 1$ and $d_i = d$ for all $2 \leq i \leq n$, then one can show that estimate (4.1) is also true if $X$ is a sufficiently nice family (Lefschetz pencil) of hypersurfaces in $\mathbf{A}^{n-1}$ parametrized by $x_1$. Again, the exact number of zeros and poles of the partial zeta function $Z_{1,d,\cdots,d}(X,T)$ seems to be out of reach due to possible cancellation of zeros. In this example, the reader may be more comfortable to use the fibred variety version and consider a nice family of smooth projective hypersurfaces. It would be important for applications to have a systematic and thorough study of the estimate of type (4.1) by fully exploiting the theory of $\ell$-adic cohomology, and by finding simple conditions on $X$ for which estimate (4.1) holds. A similar question can be asked in the more general situation when the $d_i$'s do not satisfy the dividing condition. In this case, one would need to understand the detailed geometry of the variety $Y(d_1,\cdots,d_n,X)$ in Faltings' construction and hopefully establish some strong cohomological vanishing theorems. Very little is known.

The $p$-adic variation of the partial zeta function is closely related to Dwork's conjecture [Dw2]. As an example, let us consider the special case

$$d_1 = \cdots = d_r = 1, \quad d_{r+1} = \cdots = d_n = d + p^m D,$$

where $1 \leq r < n$ is a fixed positive integer and $D$ is some positive integer. Let

$$f: X \longrightarrow \mathbf{A}^r, \ (x_1,\cdots,x_n) \longrightarrow (x_1,\cdots,x_r)$$

be the natural projection map. It is not hard to show that for certain choices of $D$, we have the congruence

$$Z_{1,\cdots,1,d+p^m D,\cdots,d+p^m D}(X,T) \equiv Z_{1,\cdots,1,d+p^{m+1}D,\cdots,d+p^{m+1}D}(X,T) \ (\mathrm{mod} \ p^m)$$

for all positive integers $m$. This $p$-adic continuity relation implies that the $p$-adic limit

$$\lim_{m\to\infty} Z_{1,\cdots,1,d+p^m D,\cdots,d+p^m D}(X,T) = L^{[d]}(0,f,T)$$

exists as a $p$-adic power series. This limiting function coincides with Dwork's slope zero L-function attached to the family $f$. In terms of $p$-adic étale cohomology, one derives the formula

$$L^{[d]}(0,f,T) = \prod_{i\geq 0} L((R^i f_! \mathbf{Z}_p)^d, T)^{(-1)^{i-1}},$$

where the L-factor on the right side denotes the L-function of the $d$-th power (iterate) of the relative $p$-adic étale cohomology $R^i f_! \mathbf{Z}_p$ with compact support. Dwork's conjecture says that the L-function $L^{[d]}(0,f,T)$ is $p$-adic meromorphic, see the survey paper [W6] for further information and [W3-5] for the proof of Dwork's conjecture.

It is our hope that the observations and questions included in this short paper would give adequate motivation for interested readers to study further properties of the partial zeta function and to explore their potential applications.

**Acknowledgement**. The author would like to thank Gerd Faltings for giving a beautiful solution to the near-rationality problem and for allowing me to include an exposition of his idea here.